\input amstex
\documentstyle{amsppt}
\document
\magnification=1200
\NoBlackBoxes
\nologo
\pageheight{18cm}

%\hfill{\it file Amywork/RMtalk.tex, draft Feb. 9, 2002}

\bigskip

\centerline{\bf VON ZAHLEN UND FIGUREN\footnotemark1}
\footnotetext{Talk at the International Conference
``G\'eom\'etrie au vingti\`eme ci\`ecle: 1930--2000'',
Paris, Institut Henri Poincar\'e, Sept. 2001. The title
is a homage to Hans Rademacher and Otto Toeplitz whose book
fascinated the author  many years ago.}

\medskip

\centerline{\bf Yu.~I.~Manin}

\medskip

\centerline{\it Max--Planck--Institut f\"ur Mathematik, Bonn}

\bigskip

\centerline{\bf \S 0. Introduction}

\medskip

{\it Geometry} is a large subfield of 
mathematics, but also
a label for  a certain mindset of a practising mathematician.
The same can be told about {\it Algebra} (understood
here broadly, as the language of mathematics,
as opposed to its content, and so including {\it Logic}.) A natural or acquired
predilection towards geometric or algebraic
thinking and respective mental objects is often expressed
in strong pronouncements, like Hermann Weyl's
exorcising ``the devil of abstract
algebra'' who allegedly struggles with
``the angel of geometry''  for the soul
of each mathematical theory. (One is reminded of an even more sweeping
truth: {\it ``L'enfer -- c'est les autres''}.) 

\smallskip

Actually, the most fascinating thing about algebra and geometry
is the way they struggle to help each other
to emerge from the chaos of non--being, from
those dark depths of subconscious where all roots of
intellectual creativity reside. What one ``sees''
geometrically must be conveyed to others in words and symbols.
If the resulting text can never be a perfect
vehicle for the private and personal vision, the vision itself
can never achieve  maturity without
being subject to the test of written speech. The latter
is, after all, the basis of the social 
existence of mathematics. 

\smallskip

A skillful use of the interpretative algebraic language possesses 
also a definite therapeutic 
quality. It allows one to fight the obsession which
often accompanies contemplation of enigmatic
Rorschach's blots of one's  imagination.

\smallskip

When a significant new unit of meaning (technically, 
a mathematical definition
or a mathematical fact) emerges from such a struggle,
the mathematical community spends some time
elaborating all conceivable implications
of this discovery. (As an example, imagine the development
of the idea of a continuous function, or a Riemannian metric,
or a structure sheaf.)
Interiorized, these implications prepare new firm ground for
further flights of imagination, and more often than not
reveal the limitations of the initial formalization
of the geometric intuition. 
Gradually the discrepancy between
the limited scope of this unit of meaning 
and our newly educated and enhanced geometric vision
becomes glaring, and the cycle repeats itself.

\smallskip

A very special role in this cyclic process is played by problems.
The importance of a new theoretical development is generally 
judged by its success (or otherwise) in throwing light on 
an old problem or two. Conversely, a problem
can stimulate the emergence of a new geometric vision
expressing a sudden perception of a hidden analogy,
as in Andr\'e Weil's famous recognition
of Lefschetz's formula in the context of algebraic
equations over finite fields.

\smallskip   

A problem/conjecture is usually represented
by a short mathematical statement allowing a yes or no answer.
``Yes'' and ``no'' do not play symmetric roles:
a positive answer to a good question 
usually validates a certain intuitive picture,
whereas a ``no'' answer often shows only the limitations of this
picture rather than its total lack of value.
A ``counterexample'' disposing of a resistant conjecture 
can have a certain sportive value, but becomes
really important only if pursued so far as to reveal some positive truth
which has been escaping our understanding for some time.

\smallskip

The geometry of XX century is a huge patchwork of
ideas, visions, problems and its solutions. A brief list of platitudes
I sketched above could be used to organize a narrative dedicated
to the contemporary history of this discipline or its separate episodes.
I have chosen instead to present at this conference a narrative based on my
current work which at its key point remains purely conjectural,
if not outright speculative.

\smallskip

I hope that the lack of a real mathematical breakthrough to report   
can be compensated in the context of this conference
by certain freshness of perception accompanying such early stages of research.
Besides, a discussion of partial successes and failures
of this enterprise can serve as an illustration
of some general issues of geometry and algebra.

\smallskip

Briefly, the research in question concerns 
explicit construction of numbers
generating abelian extensions of algebraic number fields.

\smallskip

The archetypal problem is that of understanding abelian extensions
of $\bold{Q}$. As we know after Kronecker
and Weber, the maximal abelian extension of $\bold{Q}$ is
generated by roots of unity.  Roots of unity of degree $n$ 
form in the complex plane vertices of a regular $n$--gone.
Which of these $n$--gons can be constructed using only ruler and
compass, was a famous problem solved by Gauss. His first
publication dated April 18, 1796 (and {\it Tagebuch} entry of March 30,
cf. [Ga]) is an announcement 
that a regular 17--gon has this property. Gauss was not quite 19 then; 
apparently, this discovery prompted him to dedicate his life to
mathematics.  

\smallskip
One remarkable feature
of Gauss' result is the appearance of a hidden symmetry group.
In fact, the definitions of a regular $n$--gon and  ruler
and compass constructions are given in terms
of Euclidean plane geometry and make practically ``evident''
that the relevant symmetry group is that of rigid
rotations $SO\,(2)$ (perhaps, extended by  reflections and shifts).
This conclusion turns out to be totally misleading: instead,
one should rely upon $\roman{Gal}\,(\overline{\bold{Q}}/\bold{Q})$.
The action of the latter group upon roots of unity of degree $n$ 
factors through
the maximal abelian quotient and is given by $\zeta\mapsto \zeta^k,$
with $k$ running over all $k\,\roman{mod}\,n$ with $(k,n)=1$, 
whereas the action of the rotation group is
given by $\zeta\mapsto \zeta_0\zeta$ with $\zeta_0$ running
over all $n$--th roots. Thus, $\roman{Gal}\,(\overline{\bold{Q}}/\bold{Q})$ does not conserve
angles between vertices which seem to be basic for the initial problem.
Instead, it is compatible with addition and multiplication
of complex numbers, and this property proves to be crucial.

\smallskip

This gem of classical mathematics contains in a nutshell
some basic counterpoints of my presentation:
geometry and spatial imagination vs  refined language of algebra;
physics (kinematics of solid bodies) vs number theory;
relativity of continuous and discrete. I look
at them from the modern perspective of non--commutative geometry,
inspired by several deep insights of Alain Connes:

\smallskip

(a) Connes' bold attack on the Riemann Hypothesis ([Co3], [Co4]);

\smallskip

(b) the discovery (joint with J.~Bost, [BoCo]) of a statistical
system with the symmetry group $\roman{Gal}\,(\overline{\bold{Q}}/\bold{Q})$
and spontaneous symmetry breaking at the pole $s=1$ of the
Riemann zeta;
 
\smallskip

(c) Connes' idea that approximately
finite dimensional simple central $C^*$--algebras will furnish
the missing nontrivial Brauer theory for archimedean
arithmetic infinity ([Co3], p. 72 and [Co4], p. 38).

\smallskip

I discuss some of these and related themes and add to them
the fourth suggestion:

\smallskip

(d) Real Multiplication project, in which elliptic curves
with complex multiplication are replaced by two--dimensional
quantum tori $T_{\theta}$ whose $K_0$--group (or rather, its image
wrt the normalized trace map $\bold{Z}+\bold{Z}\theta$) is a subgroup of a real quadratic field. More details will be given in a paper
in preparation [Ma4].

\bigskip 

\centerline{\bf \S 1. Real multiplication: an introduction}

\medskip 

Let $K$ be a local
or global field in the sense of number theory, that is,
a field of algebraic numbers, a field of
functions on a curve over $\bold{F}_q$, or a completion
of one of these fields. Denote by $K^{ab}$ its maximal
abelian (separable) extension. Class field theory
provides a description of the Galois group
$\roman{Gal}\,(K^{ab}/K)$ and some partial information
on its action upon $K^{ab}$. However,  
specific generators of $K^{ab}$,
together with exact action  of $\roman{Gal}\,(K^{ab}/K)$ upon them,
are not generally known. Hilbert's twelfth problem
addresses this question. Stark's conjectures, still unproved
in general, provide a partial answer: see [St1], [St2], [Ta].
Below we will define Stark's numbers directly for
$\bold{Q}$ (see 2.1) and real quadratic fields (see 4.1).

\smallskip

The most elementary, complete and  satisfactory description
is furnished by the Kronecker--Weber theorem which we already
mentioned.

\smallskip

The Kronecker--Weber machine (KW) can be successfully
imitated for at least three more classes of ground fields
$K$: complex quadratic extensions of $\bold{Q}$ (referred to as CM, for
complex multiplication theory),
finite extensions of $\bold{Q}_p$ (LT, for Lubin--Tate formal groups), and global
fields of finite characteristics (D, for Drinfeld's modules). 

\smallskip

The CM case was elaborated first, and 
all variations
on the  CM theme, including KW, can be described
according to the following scheme. 
(Only the Lubin--Tate case slightly diverges from the general
pattern at this point, see below).

\smallskip

Let $K$ be a global field as above,
$O\subset K$ its appropriate subring with quotient field $K$. 
We denote by $\bold{A}$
the analytic (or formal) additive group defined over $K$
and construct the quotient group $\bold{A}/O$ where $O$
is embedded as a subgroup of $K$--points of $\bold{A}$.
By construction, multiplication by $O$ induces on $\bold{A}/O$ endomorphisms
in an appropriate category of $K$--groups: this is
the {\it geometric} multiplication.

\smallskip

We construct a function $t$ on $\bold{A}/O$ taking
algebraic (over $K$) values at the $O$--torsion points
$\xi \in K/O\subset \bold{A}/O$. The action of a Galois
group  upon these torsion
points will have an abelian image, if it is shown to
commute with the geometric action of $O$. This presupposes
a careful study of the fields of definition of various objects involved. 

\smallskip

Class field
theory plus special properties of the function
$t$ allows us then to establish that the values of
$t$ on $K/O$ generate (almost) all of $K^{ab}/K.$ 
Here are some more details.

\smallskip 
 
{\bf 1.1. Case (KW).} Here $K:=\bold{Q}$,
$O:=\bold{Z}$. For $\bold{A}$ we will take the additive
group considered as a complex--analytic space supplied
by an additive coordinate $z$. The
map $a:\,z\mapsto e^{2\pi iz}$ identifies $\bold{A}/O=\bold{G}_a^{an}/\bold{Z}$
with the (analytic) multiplicative group
$\bold{G}_m^{an}.$ The latter, of course, has a canonical
algebraic structure, with a multiplicative
coordinate which we will denote $1+t$, $a^*(1+t)=e^{2\pi iz}$,
and this $t$ will be the function mentioned in the general
description. The lattice $\bold{Z}$ acts upon
$\bold{G}_m$ by $1+t\mapsto (1+t)^n$: this is the geometric
action. The values of $1+t$ 
at the torsion points are all roots of unity.
The Galois group commutes with geometric semigroup and
therefore can be identified with its profinite completion
$\widehat{\bold{Z}}^*$.

\smallskip

{\bf 1.2. Case (CM).} Here we denote
the ground field (former $K$) by $L$ and assume that it is an imaginary
quadratic extension of $\bold{Q}.$ For $O$ we take the
ring $O_L$ of all integers in $L$, and for $\bold{A}$
the complex analytic additive group as above.
The homomorphism of analytic groups $a:\, \bold{G}_a\to X_L$ 
with the kernel $O_L$ is the universal covering of the
elliptic curve $X_L$. Over $\bold{C}$, its field of
algebraic functions is generated by the Weierstrass
function $\wp (z, O_L)$ and its derivative.
The action of $O_L$ upon $X_L$,
in particular,  produces rational expressions
for $\wp (\alpha z, O_L)$, $\alpha\in O_L,$
in terms of $\wp (z, O_L)$, $\wp^{\prime} (z, O_L)$.

\smallskip

{\it (i) First description of $L^{ab}$.}
An important difference from the (KW)--case
is that elliptic curves, unlike $\bold{G}_M$,
are not rigid but occur in families.
In particular, the isomorphism class of $X_L$  over $\bold{C}$
is determined by the invariant
$j(O_L)\in\bold{C}$. It turns out that $L(j(O_L))$
is exactly the maximal unramified abelian extension
of $L$, i.e. its Hilbert class field.

\smallskip

Over $L(j(O_L))$, the maximal abelian extension
of $L$ can be generated by the values of the function $t$
such that $a^*(t):=\wp ^{-u}(z, O_L)$, where
$2u$ is the order of the group $O_L^*$ of units in $O_L$.
Generally, $u=1$; but $u=2$ (resp. $u=3$) for 
$L=\bold{Q}(e^{\pi i/2})$ (resp. $L=\bold{Q}(e^{\pi i/3})$).

\smallskip

A more geometric description of $t$ is this:
$X_L/O_L^*$ is a projective line, and $t$
is a coordinate on it vanishing at the image of
the zero point of $X_L$.

\smallskip

The appearance of $X_L/O_L^*$ rather than $X_L$
itself in this description might seem a minor matter.
However, in the conjectural picture of
real multiplication it will acquire a great importance,
because $O_L^*$ for {\it real} quadratic fields
is always infinite, and the study of the respective quotient
will again require tools of noncommutative geometry.

\smallskip

{\it (ii) Second description of $L^{ab}$.} Instead
of using only the curve $X_L$ with complex multiplication
by the whole $O_L$, one can consider the family
of all elliptic curves $Y_L$ whose endomorphism
ring is an order in $O_L$. It turns out
that their absolute invariants $j(Y_L)$
generate a big subfield of $L^{ab}$: to obtain
all of $L^{ab}$, it remains however to make
an additional 2--extension (of infinite degree).
The curves $Y_L$ are all isogenous to $X_L$,
and since the points of finite order of $X_L$
``essentially'' generate $L^{ab}$, this gives 
an intuitive explanation of why this might be so.

\smallskip

However, the two description differ not only by the
details of the calculations. What is important in the
second picture, is the implicit appearance 
of the tower or stack of modular curves,
parametrizing all elliptic curves, rather than
of only one elliptic curve $X_L$.

\smallskip

This tower consists of the compactified algebraic models
of analytic spaces $H/\Gamma$, $\Gamma$, where
$H$ is the upper complex half--plane, and $\Gamma$
runs over congruence subgroups of $PSL(2,\bold{Z}).$

\smallskip

For a readable review of two approaches, see [Se] and [Ste].

\smallskip

{\bf 1.3. Case (D).} Here the ground field
is the field of functions $k=\bold{F}_q(C)$ on a smooth algebraic
curve defined over the finite base field with $q$  elements.
We choose a point $\infty \in C(\bold{F}_q)$  as a part
of the structure and denote by $O_k$ the ring
of regular functions on the affine curve $C\setminus\{\infty\}.$
For $\bold{A}$ we take the additive group considered as
an analytic group over the quotient field $k_{\infty}$ 
of the completion of $\Cal{O}_{C,\infty}.$

\smallskip

In this situation we have an analogue of exponential function
on $\bold{A}$:
$$
e(z):=z\prod_{\alpha\in O_L\setminus\{0\}} (1-\frac{z}{\alpha})
$$ 
 This function is entire. It is however additive
rather than multiplicative: $e(z_1+z_2)=e(z_1)+e(z_2).$ 
A large abelian extension of $k$, as in the
KW--case, is generated by the values of $e(\lambda z)$ at the
``points of finite $O_k$--order'', where $\lambda$
is an appropriate analogue of $2\pi i.$ The action of
$O_k$ on $\bold{G}_a$
can be alternatively described by the
embedding $O_k$ into the (algebraic) crossed product
of $k_{\infty}$ with the semigroup of non--negative powers of
the Frobenius endomorphism $x\mapsto x^p$ acting upon $k_{\infty}$.

\smallskip 

This construction can be considerably generalized by
replacing $\lambda O_k$ in the construction of the
exponential function with any $O_k$ lattice
of rank $d\ge 1.$ We get in this way
the notion
of the Drinfeld module of rank $d$. The case $d=1$ 
produces abelian extensions
and is similar both to (KW) and (CM) cases.
The rank 2 case looks like a version of the theory
of elliptic curves, and an appropriate modification
of the general  $d$ case 
produces the Langlands type
description of the algebraic closure of $K$.

\smallskip

A thorough historical and mathematical discussion of the
theories (KW), (CM) and (D) can be found in the book [V].

\smallskip

{\bf 1.4. Case (LT).} Here $K$ is a finite extension of $\bold{Q}_p$,
$O=O_K$ is the ring of integers in $K$, and the group $G_{LT}$
denoted formerly by $\bold{A}/O$ above, is an one--dimensional formal group
over $O_K$ equipped with a homomorphism $O_K\to \roman{End}\,G_{LT}$
sending any $a\in O_K$ to a formal map with the linear term $x\mapsto ax$.

\smallskip

I do not know whether $G_{LT}$ can be interpreted as a quotient
$\bold{A}/O$ in an appropriate category, although the theory of logarithmic
functions developed in the context of
$p$--adic Hodge--Tate theory indicates that such an interpretation
might exist. Lubin and Tate construct $G_{LT}$ directly
using an ingenious calculation with formal series
providing simultaneously a description of the category
of such formal groups.

\smallskip

{\bf 1.5. Real multiplication (RM) and noncommutative geometry.}
The simplest class of fields $M$ for which no direct description
of $M^{ab}$ is known consists of
real quadratic extensions of $\bold{Q}.$ 

\smallskip

Below we sketch a possible approach to this problem
via non--commutative geometry and suggest its
relation with the Stark conjectures. The idea is
straightforward enough. An elliptic curve with
complex multiplication by, say, the maximal order
$O_L$ in a complex quadratic field $L$
has a complex analytic model $\bold{C}/O_L$,
where $O_L$ is considered as lattice in $\bold{C}.$
Similarly, one could try to imagine a space
$\bold{R}/O_M$ where $O_M$ is the ring of integers
of a real quadratic field $M$. Of course,
$O_M$ is not discrete anymore, but the first principle
of the non--commutative philosophy says that
one should not shy away from such situations:
$\bold{R}/O_M$ exists as a non--commutative space
represented e.g. by a two--dimensional torus
$T_{\theta}$ where $O_M=\bold{Z}+\bold{Z}\,\theta .$
Recall that the $C^*$--algebra of $T_{\theta}$
is the universal algebra generated by two unitaries
$U,V$ with the commutation rule $UV=e^{\pi i\theta}VU.$
This algebra is the crossed product of $C(\bold{R}/\bold{Z})$ and
the irrational rotation automorphism of  $\bold{R}/\bold{Z}$ induces by
the shift $t\mapsto t+\theta$. Crossed products generally
serve to represent quotients with respect to
``bad'' (and occasionally good) equivalence relations:
cf, [Co1], pp. 85--91 for a very livrly account of this principle.
(With a hindsight, we noted  that in the construction
of Drinfeld modules the central role was played
by a universal crossed product of an algebra of power
series and the semigroup generated by the Frobenius endomorphism.)

\smallskip

As the next step, we want to make sense of the statement
that $T_{\theta}$ admits real multiplication by $O_M$.
This is obvious enough for $\bold{R}/O_M$ understood
set--theoretically: the action of $O_M$
is simply induced by the multiplication $O_M\times\bold{R}\to \bold{R}$. 
But crossed product $C^*$--algebras are not functorial in
any naive sense. One way to deal with this difficulty is to
replace the usual homomorphisms of associative rings by
Morita functors between the categories of their (say, right) modules.
Morita functors are given by bimodules, and  composition of functors
corresponds to the tensor product of bimodules. 

\smallskip

Similarly, a coarse moduli space of two--dimensional quantum tori
up to Morita equivalence can be seen as a
quotient $PGL (2,\bold{Z})\setminus\bold{P}^1(\bold{R})$; various rigidities
lead to its modular covers. Basics of topology and function theory
of these spaces were studied in [MaMar].

\smallskip

Finally, a key issue for arithmetics is the question of
fields of definition of various objects:
of the ``non--commutative elliptic curves with real 
multiplication $T_{\theta}$, ''  of the Galois--Hopf automorphism algebras
of such objects (replacing their points of finite order) etc. 
An implication is that we need algebraic geometric,
finitely or countably generated objects, preferably over
$\bold{Z}$, from which the functional
analytic structures like $C^*$--algebras can be obtained by extension
of the base field and an appropriate completion. 

\smallskip

In the next section I will address some of these challenges. 
I will start it with a discussion of several constructions
in which algebraic numbers, $C^*$--algebras and related
noncommutative rings appear  in combinations that are likely
to be susceptible to generalizations.

\bigskip

\centerline{\bf \S 2. Noncommutative geometry and arithmetic}

\medskip

Oversimplifying, one can say
that in commutative geometry algebraic numbers
appear as values of algebraic functions, whereas
in noncommutative geometry they appear as values of 
traces of projections, or more generally values
of appropriate states on observables. In both cases,
a control of the action of the Galois group
is gained, if this action commutes with an action
of certain ``geometric'' endomorphisms, or correspondences,
whenever the latter are defined over the ground field.

\smallskip

We will deal with three basic situations:

\smallskip

(A) V.~Jones theory of indices of subfactors
(cf. [Jo1], [Jo2], [GoHaJo]) and its further developments.

\smallskip

(B) Bost--Connes ``spontaneous symmetry breaking'', with
the symmetry group $\roman{Gal}\, (\overline{\bold{Q}}/\bold{Q})$
(cf. [BoCo]). 
\smallskip

(C) V.~Drinfeld's embedding of $\roman{Gal}\,(\overline{\bold{Q}}/\bold{Q})$
into the Grothendieck--Teichm\"uller group: cf. [Dr], [DE], [LoS].

\medskip

{\bf 2.1. Jones indices.} Consider first the Temperley--Lieb algebra
$TL\,(n+1,\tau )$ which arose in statistical physics. Over 
the central subring $\bold{Z}[\tau ]$,
it is generated by the idempotents $e_1,\dots ,e_n$
satisfying the relations $e_ie_j=e_je_i$ for
$|i-j|\ge 2$ and $e_ie_{i\pm 1}e_i=\tau e_i.$
It has a finite dimension, and is semisimple for
generic $\tau$.
For a critical
value of the form $\tau^{-1}=4\,\roman{cos}^2\frac{m\pi}{n+1}$, 
$TL\,(n+1,\tau )$ fails to be semisimple.

\smallskip

It is interesting that Stark's numbers for the ground field $\bold{Q}$
([Ta], p.~79) 
 differ from these critical numbers only
by a sign change and shift by 4:
$$
\roman{exp}\,(-2\zeta_{(m,n)}^{\prime}(0))=4\,\roman{sin}^2\frac{m\pi}{n}
$$
where
$$
\zeta_{(m,n)}(s) :=\sum_{k\in m+n\bold{Z}}|k|^{-s}.
$$
In particular, TL--critical values and Stark numbers
generate the maximal real subfield of $\bold{Q}^{ab}$ and are stable
with respect to the Galois group of $\overline{\bold{Q}}.$

\smallskip

One can extend the base ring of TL--algebras to $\bold{C}$, define the 
*--involution on them by $e_i^*=e_i$ (so that $e_i$ become projections),
and then pass to the inductive limit 
with respect to the obvious injections $TL\,(n+1,\tau )\mapsto
TL\,(n+2,\tau )$. It turns out that the resulting
algebra admits an involutive representation in a complex Hilbert space
iff either $\tau^{-1}=4\,\roman{cos}^2\frac{\pi}{n}$ for some $n\ge 3$,
or $\tau^{-1}\ge 4.$

\smallskip

This statement constitutes an essential part of the
famous result due to V.~Jones who proved that
for any pair of $\roman{II}_1$--factors $N\subset M$
the index $[M:N]$ lies in the set $\{4\,\roman{cos}^2\frac{\pi}{n},\,n\ge 3\}\,
\cup\,[4,\infty ]$,  and that for the hyperfinite $M$ this is exactly the set of values of indices. Index itself can be defined as a value of the Hattori--Stallings 
rank, or as a measure of the growth rate of the minimal
number of generators of the left $N$--module 
$M^{\otimes n}:=M\otimes_N\dots\otimes_NM$ as $n\to\infty$.

\smallskip

Notice that although the discrete part of the values of the index
constitutes a part of the TL--spectrum, it is not 
$\roman{Gal}\,(\bold{Q}^{ab}/\bold{Q})$--invariant.
The reason is that by passing to the $C^*$--limit, we have 
implicitly chosen a non--archimedean  valuation of
$\bold{Q}^{ab}$.

\medskip

{\bf 2.2. Bost--Connes spontaneous symmetry breaking.} In the remarkable
paper [BoCo] the action of $\roman{Gal}\,(\bold{Q}^{ab}/\bold{Q})$
on roots of unity appears in yet another setting. Instead of (the union of)
the Temperley--Lieb algebras, consider the Hecke algebra $\Cal{H}$
with involution  over 
$\bold{Q}$
given by the following presentation. The generators
are denoted $\mu_n,\,n\in \bold{Z}_+$, and $e(\gamma )$, $\gamma\in
\bold{Q}/\bold{Z}.$ The relations are
$$
\mu_n^*\mu_n=1,\quad \mu_{mn}=\mu_m\mu_n,\quad
\mu_m^*\mu_n=\mu_n\mu_m^*\ \roman{for}\ (m,n)=1;
$$
$$
e(\gamma)^*=e(-\gamma ),\quad e(\gamma_1+\gamma_2) =
e(\gamma_1)\,e(\gamma_2);
$$
$$
e(\gamma )\mu_n= \mu_n e(n\gamma ),\quad 
\mu_n e(\gamma )\mu_n^*=\frac{1}{n}
\sum_{n\delta =\gamma} e(\delta ).
$$
The id\`ele class group $\widehat{\bold{Z}}^*$ of $\bold{Q}$ acts upon $\Cal{H}$
in a very explicit and simple way: on $e(\gamma )$'s
the action is induced by the multiplication $\widehat{\bold{Z}}^*\times \bold{Q}/\bold{Z}
\to \bold{Q}/\bold{Z}$, whereas on $\mu_n$'s it is identical.

\smallskip

The algebra $\Cal{H}$ admits an involutive representation $\rho$ in $l^2
(\bold{Z}_+)$: denoting by $\{\epsilon_k\}$ the standard basis
of this space, we have
$$
\rho (\mu_n)\,\epsilon_k = \epsilon_{nk},\quad
\rho (e(\gamma ))\epsilon_k= e^{2\pi ik\gamma}\epsilon_k.
$$
From this, one can produce the whole
$\roman{Gal}\,(\bold{Q}^{ab}/\bold{Q})$--orbit $\{\rho_g\}$ of 
such representations,
applying $g\in \roman{Gal}\,(\bold{Q}^{ab}/\bold{Q})$
to all roots of unity occuring at the right hand sides of
the expressions for  $\rho (e(\gamma ))\epsilon_k$.
All these representations can be canonically extended
to the $C^*$--algebra completion $C$ of $\Cal{H}$
constructed from the regular representation of $\Cal{H}$.
Let us denote them by the same symbol $\rho_g$.

\smallskip

To formulate the main theorem of [BoCo], we need some more explanations.
The algebra $C$ admits a  canonical action of $\bold{R}$,
which can be interpreted as {\it time evolution} 
represented on the algebra of observables. This is
a general (and deep) fact in the theory of $C^*$--algebras,
but for $C$ the action of $\bold{R}$ can be quite explicitly
described on the generators. Let us denote by $\sigma_t$
the action of $t\in\bold{R}.$ {\it A KMS${}_{\beta}$ state
at inverse temperature} $\beta$ on $(C,\sigma_t)$
is defined as a state $\varphi$ on $C$ such that for any
$x,y\in C$ there exists a bounded holomorphic function
$F_{x,y}(z)$ defined in the strip $0\le \roman{Im}\,z\le
\beta$ and continuous on the boundary, satisfying
$$
\varphi (x\sigma_t(y)) = F_{x,y}(t),\ 
\varphi (\sigma_t(y)x) = F_{x,y}(t+i\beta ).
$$
Now denote by $H$ the positive operator on
$l^2(\bold{Z}_+)$: $H\epsilon_k = (\roman{log}\,k)\,\epsilon_k.$
Then for any $\beta >1$, $g\in \roman{Gal}\,(\bold{Q}^{ab}/\bold{Q})$
one can define a KMS${}_{\beta}$ state $\varphi_{\beta,g}$
on $(C,\sigma_t)$ by the following formula:
$$
\varphi_{\beta,g}(x):=\zeta (\beta )^{-1}\,\roman{Trace}\,
(\rho_g(x)\,e^{-\beta H}),\quad x\in C
$$
where $\zeta$ is the Riemann zeta--function.
The map $g\mapsto \varphi_{\beta,g}(x)$ is a homeomorphism
of $\roman{Gal}\,(\bold{Q}^{ab}/\bold{Q})$ with the space of extreme
points of the Choquet simplex of all KMS${}_{\beta}$ states.

\smallskip

To the contrary, for $\beta <1$ there is a unique KMS${}_{\beta}$ state.
This is a remarkable ``arithmetical symmetry breaking'' phenomenon.

\medskip

{\bf 2.3. $\roman{Gal}\,(\overline{\bold{Q}}/\bold{Q})$ and 
the Grothendieck--Teichm\"uller group.} Here we briefly describe the setting
studied in [Dr]. Consider the following tower of fields
$$
\bold{Q}(t)\subset \overline{\bold{Q}}(t)\subset F
$$ 
where $F$ is the maximal algebraic extension of $\overline{\bold{Q}}(t)$
ramified only at $t=0,1,\infty .$ The Galois
group $\roman{Gal}\,(F/\overline{\bold{Q}}(t))$ is an extension
of $\Cal{G}:=\roman{Gal}\,(\overline{\bold{Q}}/\bold{Q})$ by
$\Pi:=\roman{Gal}\,(F/\overline{\bold{Q}}(t)).$
Hence $\Cal{G}$ acts upon $\Pi$ by outer automorphisms.
Now, $\Pi$ is a profinite completion of the fundamental
group  $\pi_1(\bold{P}^1(\bold{C})\setminus\{0,1,\infty \})$.
Drinfeld prefers to work with the ``formal'' monodromy group
of the differential equation
$$
G^{\prime}(z)=\frac{1}{2\pi i}\,\left(\frac{A}{x}+\frac{B}{z-1}\right)\,
G(z)
$$
where $A,B$ are non--commuting symbols, and with Knizhnik--Zamolodchikov
(KZ) scattering data for such an equation treated
in a purely algebraic way. This allows him to connect
the Galois action with his study of quasi--Hopf algebras.

\smallskip

\medskip

{\bf 2.4. Summary.} The three constructions briefly summarized
above at the moment do not form a part of a coherent
picture. However, they shed some light upon each other.

\smallskip

The Temperley--Lieb algebras were recently understood
as a kind of Galois symmetry objects (quantum groupoids)
responsible
for the classification of pairs of factors
$N\subset M$ with finite width and index:
see the review [NiVa] and other references quoted therein.
Apparently, the action of the cyclotomic
Galois group on them comes from a setting
similar to that described in [Dr], which however
has to be made explicit yet. Although in the
Grothendieck--Drinfeld picture the stress was usually made
on the dessins d'enfant rather than unwieldy
non--commutative objects, this might be only one
of several possibilities.

\smallskip

Approximately finite dimensional algebras are 
completed inductive limits of split central
semisimple algebras, i.e. direct sums of matrix algebras.
It would be interesting to develop the theory
of inductive limits of general central 
semisimple algebras over non--closed, in particular,
number fields.

\smallskip

No explicit connection was made between Bost--Connes
theory and Jones theory.  Are there any? Since zeta function
appears explicitly in [BoCo] as a partition function, cyclotomic Stark numbers
might surface naturally in  this model.

\smallskip

Finally, all three theories are absolute in the sense
that their ground field is $\bold{Q}$. It is striking
that their relative analogs (over arbitrary number fields $K$)
are not known.

\smallskip 

To be more precise, some rather straightforward 
extensions of [BoCo] were studied in
[HaL], [ArLR], [Coh1], [Coh2] but none of them has the $\pi_0$ of 
the id\`ele class group of $K$
as a symmetry group  when $K\ne \bold{Q}.$
In fact,
the archimedean part of the id\`eles is not properly accommodated,
which returns us to the problem of the right Brauer theory
for archimedean primes.

\smallskip 

One can expect that a real understanding
of the classical abelian class field theory
will  depend on such a generalization.
 
\bigskip

\centerline{\bf \S 3. Real multiplication of quantum tori: geometry}

\medskip

{\bf 3.1. Lattices and pseudolattices.} The category of elliptic curves
over $\bold{C}$ is equivalent to the category of
period lattices. Formally, a lattice
is a discrete and cocompact embedding $j:\,\Lambda \to V$ where $\Lambda$
is isomorphic to $\bold{Z}^2$ and $V$ is an one-dimensional
complex vector space; morphisms are linear maps $V\to V^{\prime}$
sending $j(\Lambda )$ into $j^{\prime}(\Lambda^{\prime} )$.
The functor establishing equivalence sends $(\Lambda ,V, j)$
to the complex torus $V/j(\Lambda ).$ 
\smallskip

Let us denote by $\Lambda_{\tau}$ the lattice $\bold{Z}+\bold{Z}\tau ,$
$\tau\in \bold{C}\setminus\bold{R}.$
Then any lattice is isomorphic to some
$\Lambda_{\tau}$, and each non--zero morphism $\Lambda_{\tau^{\prime}} \to \Lambda_{\tau}$ is represented
by a non--degenerate matrix 
$$
g=\left(\matrix a&b\\c&d\endmatrix\right)\in M\,(2,\bold{Z})
$$
such that
$$
\tau^{\prime} =\frac{a\tau +b}{c\tau +d},\quad
$$ 
It  is an isomorphism, iff $g\in GL(2,\bold{Z})$.
Thus the coarse moduli
space classifying lattices up to an isomorphism is
$$
PGL(2,\bold{Z})\setminus (\bold{P}^1(\bold{C})\setminus \bold{P}^1(\bold{R}))
$$
which is usually written as $PSL(2,\bold{Z})\setminus H$ where $H$
is the upper half--plane. 

\smallskip

The map $j$ can degenerate so that $j$ acquires a kernel of rank one;
the image remains discrete but is not cocompact.
The coarse moduli space admitting such degenerations is $PSL(2,\bold{Z})\setminus (H\cup \bold{P}^1(\bold{Q}))$: a cusp is added. This cusp (or several cusps if
a rigidity is imposed) can be seen on algebraic geometric
models of the modular curves and canonical elliptic fibrations over them:
they essentially correspond to the degenerations of
elliptic curves with (stably) multiplicative connected component
of the closed fiber.
\smallskip

Classes of the irrational
real points constitute another part of the boundary 
$PGL(2,\bold{Z})\setminus \bold{P}^1(\bold{R})$.
This part classifies degenerations
of lattices which I will call {\it pseudolattices}:
a pseudolattice is an embedding $j:\,L \to V$
where $L$ is isomorphic to $\bold{Z}^2$, $V$ is
an one--dimensional $\bold{C}$--space as above,
and the closure of $j(V)$ is a real line.
Morphisms are defined exactly as for lattices.
Denoting by $L_{\theta}$ the pseudolattice $\bold{Z}+\bold{Z}\theta$,
we see as above that any pseudolattice
is isomorphic to an  $L_{\theta}$, and morphisms correspond
to the fractional linear transformations of $\theta$'s 
with integral coefficients.

\smallskip

These degenerations are invisible in algebraic geometry
because $V/j(L)$ makes no sense as an algebraic
or analytic curve. But the machinery of
noncommutative geometry of Connes is designed
to deal with such spaces. Choosing $L_{\theta}$ as a representative of the
respective isomorphism class, we can naively replace
$\bold{C}/(\bold{Z}+\bold{Z}\theta )$ by
$\bold{C}^*/(e^{2\pi i\theta})$ (``Jacobi uniformization''),
and then interpret the last quotient
as an ``irrational rotation algebra'', or
two--dimensional quantum torus $T_{\theta}.$ We recall
that this torus
is (represented by) the universal $C^*$--algebra $A_{\theta}$
generated by two unitaries $U,V$ with the commutation rule
$UV=e^{2\pi i\theta}VU.$ 

\smallskip

The next task is to define morphisms between these quantum tori,
so that we could compare their category with that
of pseudolattices. Already isomorphisms present a problem:
we want invertible fractional linear transforms to produce
isomorphic quantum tori. M.~Rieffel's seminal
discovery (cf. [Ri3]--[Ri5], [RiSch]) was that to this end we should consider
Morita equivalences between appropriate categories of modules 
as isomorphisms between the tori themselves.
Taking this lead, we will formally introduce the general
Morita morphisms of associative rings, stressing those
traits of the formalism that play a central role in the
structure theory of quantum tori (of arbitrary dimension).

\smallskip

For a considerably more sophisticated version of such a theory
for von Neumann algebras, see [Co1], VB, and a report
[Jo3]. Physical motivations can be found
in [Sch1], [Wi], [CoDSch].

\medskip

{\bf 3.2. Morita category and projective modules.}  Let $A,B$ be two associative rings.
A {\it Morita morphism} $A\to B$
by definition, is the isomorphism class of a bimodule 
${}_AM_B$, which is projective and finitely generated
separately as module over $A$ and  $B$.

\smallskip

The composition of morphisms is given by the tensor product
${}_AM_B\otimes {}_BM^{\prime}_C,$ or ${}_AM\otimes {}_BM^{\prime}_C$
for short.

\smallskip  

If we associate to ${}_AM_B$
the functor
$$
\roman{Mod}_A\to \roman{Mod}_B:\ 
N_A\mapsto N\otimes_AM_B,
$$
the composition of functors will be given by the tensor product,
and isomorphisms of functors will correspond to the isomorphisms
of bimodules.

\smallskip

We imagine an object $A$ of the (opposite) Morita category  as a  
noncommutative space, right $A$--modules as sheaves on this
space, and the tensor multiplication by ${}_AM_B$
as the pull--back functor.
\smallskip

Two bimodules ${}_AM_B$ and ${}_BN_A$ supplied with two
bimodule isomorphisms ${}_AM\otimes_BN_A \to {}_AA_A$ 
and ${}_BN\otimes_AM_B\to {}_BB_B$ define
mutually inverse Morita isomorphisms (equivalences) between
$A$ and $B$. The basic example
of this kind is furnished by $B=\roman{Mat}\,(n,A)$,
$M={}_AA^n{}_B$ and $N={}_BA^n{}_A$.

\smallskip

Projective right $A$--modules up to isomorphism
are exactly ranges of idempotents in various
matrix rings $\roman{Mat}\,(n,A)$ acting from the left upon
(column) vector modules $A^n$.

\smallskip

We prefer to work with all $n$ simultaneously.
So we will denote by $\Cal{M}A$ the ring of infinite matrices
$(a_{ij})\, i,j\ge 1,\,a_{ij}\in A,$ $a_{ij}=0$ for $i+j$
big enough (depending on the matrix in question). Similarly, denote
by $A^{\infty}$ the left $\Cal{M}A$--module of infinite columns $(a_i),\, i\ge 1,$ with coordinates in $A$, almost all zero.

\smallskip

Denote by $Pr_A$ the category of finitely generated right $A$--modules.
Denote by $pr_A$ the category whose objects are
projectors $p\in \Cal{M}A,\,p^2=p$, whereas morphisms are
defined by
$$
\roman{Hom}\,(p,q):=q\,\Cal{M}A\,p
$$
with the composition induced by the multiplication
in $\Cal{M}A$.

\smallskip

There is a natural functor $pr_A\to Pr_A$ defined on objects
by $p\mapsto pA^{\infty}$. In order to define it on
morphisms, we remark that morphisms $pA^{\infty}\to qA^{\infty}$
can be naturally described by matrices in the following way.
Clearly, $p A^{\infty}$ contains the columns
$p_k$ of $p$ which generate $p A^{\infty}$ as a right $A$--module.
We can then apply any $\varphi :\, p A^{\infty}\to q A^{\infty}$
to all $p_k$ and arrange the resulting vectors into a matrix
$\Phi\in \Cal{M}A$ with $k$--th column $\varphi (p_k).$
One checks that $\Phi p=\Phi$, and since also
$q \Phi =\Phi$, we have $\Phi\in q\Cal{M}Ap$. Conversely,
any such matrix determines a unique morphism $pA^{\infty}\to qA^{\infty}$.

\medskip

\proclaim{\quad 3.2.1. Claim} (a) The functor
$pr_A\to Pr_A$ described above is an equivalence of categories.

\smallskip

(b) If $p A^{\infty}$ is isomorphic to $q A^{\infty}$, then
$p-q\in [\Cal{M}A,\Cal{M}A].$
\endproclaim 

\smallskip

{\it A trace} of $A$ is any homomorphism of additive groups $t:\,A\to G$
vanishing on commutators; by definition, it
factors through the universal trace $A\to A/[A,A]$.
Combining it with the matrix trace, we get
its canonical extension to $\Cal{M}A$. From the Claim above
it follows that $t(p)$ depends only on the
isomorphism class of $pA^{\infty}.$ The class $p\,\roman{mod}\, [A,A]$
is called {\it the Hattori--Stallings rank} of $p A^{\infty}$.

\smallskip

We define $K_0(A)$ as the Grothendieck group
of $Pr_A$. If $N_A\in Pr_A$,
$[N_A]$ denotes its class in $K_0(A).$
If $N_A$ is the range of an idempotent $p$
and $t$ is a trace, $t(p)$ depends only on $[N_A]$ and is
additive on exact triples, hence 
$t$ becomes a homomorphism of $K_0(A)$
(it is called dimension in the theory of von Neumann algebras).

\smallskip

Assume now that
$A$ is endowed with an additive
(linear or antilinear in the case of algebras) involution $a\mapsto a^*$, $(ab)^*=b^*a^*,$
$a^{**}=a.$
It extends to matrix algebras: $(B^*)_{ij}:=B_{ji}^*.$
similarly, it extends to $A^{\infty}\to A_{\infty}$
and $ A_{\infty}\mapsto A^{\infty}$, compatibly
with the module structures.

\smallskip

In such a context, it makes sense to consider only those
projective modules that are ranges of {\it projections},
that is, $*$--invariant idempotents. 

\medskip

{\bf 3.3. Morita category of two--dimensional quantum tori.}
By definition, two--dimensional quantum tori are objects
of the category $\Cal{Q}\Cal{T}$ whose morphisms are
isomorphism classes of bimodules ${}_AM_B$
corresponding to $*$--invariant projections. The algebra
$A_{\theta}$
has a unique trace $t_A$ which is normalized
by the condition $t_A(1)=1$. We can now define
a functor $K$ from $\Cal{Q}\Cal{T}$ to the category
of pseudolattices $\Cal{P}\Cal{L}$. 
On objects, we put:
$$
K(T)=(L_A,V_A,j_A,s_A).
$$
Here
$L_A:=K_0(A)$, the $K_0$--group
of the category of right projective $A$--modules
(as above, given by projections in finite matrix algebras over $A$); 
$V_A$ is the target group of the universal trace on $A$,
that is, the quotient space of $A$ modulo the completed
commutator subspace $[A,A]$. Furthermore, 
$j_A=t_A:\,K_0(A)
\to V_A$
is this universal trace extended to matrix algebras; its value
on the class of a module, as we already explained, is its value
at the respective projection. The pseudolattice $K(A_{\theta})$
comes in fact with an additional structure
which we will call {\it orientation}: namely, the cone
of effective elements in $K_0(A_{\theta})$. The information
it carries is exactly the choice of an orientation
of the real closure of $j_A(L_A).$
 
\smallskip
On morphisms, we put
$$
K({}_AM_B)([N_A]):=[N\otimes_AM_B].
$$
It takes some work to show that this functor is well defined.
Clearly, it is compatible with orientations.

\smallskip

Unlike the case of elliptic curves, $K$ is not quite an equivalence
of categories. It is essentially
surjective
on objects and morphisms of oriented pseudolattices
conserving orientations. However, it glues together
some Morita morphisms. The most clear cut
case of this is furnished by Morita
equivalences: if $\otimes {}_AM_B$
and $\otimes {}_AM_B^{\prime}$ produce Morita equivalences,
and induce the same isomorphism of pseudolattices,
these functors differ by an automorphism
of the category $\roman{Mod}_A$ which is induced by
an automorphism of the ring $A$.

\smallskip

Thus, in order to achieve an equivalence
of categories, one should count as equivalent many
morphisms of noncommutative tori induced by ring homomorphisms,
contrary to the intuition educated on affine schemes.
In fact, when the trace is unique,
all the standard automorphisms of $A_{\theta}$ act trivially
on $K_0(A_{\theta})$. But in any case, all endomorphisms
of pseudolattices keeping orientation
are induced by bimoduli, which is a basis of
real multiplication to which we finally turn.
  
\medskip

{\bf 3.4. Complex Multiplication and Real Multiplication.}
Endomorphisms
of a lattice $\Lambda$ and of the respective elliptic curve form a ring
which always contains $\bold{Z}.$ The situations
when it is strictly larger than $\bold{Z}$ can be described
as follows.

\smallskip

(a) $\roman{End}\,\Lambda \ne \bold{Z}$
iff there exists a complex quadratic subfield $K$ of $\bold{C}$ 
such that $\Lambda$ is isomorphic to a lattice contained in $K.$

\smallskip

(b) If this is the case, denote by $O_K$ the ring of integers of
$K$. There exists a unique integer $f\ge 1$ (conductor) such that
$\roman{End}\,\Lambda =\bold{Z}+fO_K=:R_f,$ and $\Lambda$ is
a projective module of rank 1 over $R_f$. Every $K$, $f$
and a projective module over $R_f$ come from a lattice.

\smallskip

(c) If lattices $\Lambda$ and $\Lambda^{\prime}$ have the same
$K$ and $f$, they are isomorphic if and only if
their classes in the Picard group $\roman{Pic}\,R_f$
coincide.

\smallskip

Automorphisms of a lattice generally form a group $\bold{Z_2}$
($\psi$ is multiplication by $\pm 1$.) However, integers
of two imaginary quadratic fields obtained by adjoining
to $\bold{Q}$ a primitive root of unity of degree 4 (resp. 6)
furnish examples of lattices with automorphism group
of order 4 (resp. 6). Only these two fields produce
lattices with such extra symmetries.

\medskip

Pseudolattices with  real multiplication admit
a completely similar description.

\smallskip

Endomorphisms of a pseudolattice $L$
(we omit other structures in notation if there is no
danger of confusion) form a ring
$\roman{End}\,L$
with 
composition as multiplication.  It contains $\bold{Z}$
and comes together with its embedding in $\bold{R}$
as $\{ a\in\bold{R}\,|\, a j(L)\subset j(L) \}.$

\smallskip

(a${}^{\prime}$) $\roman{End}\,L\ne \bold{Z}$
iff there exists a real quadratic subfield $K$ of $\bold{R}$ 
such that $L$ is isomorphic to a pseudolattice contained in $K.$

\smallskip

(b${}^{\prime}$) If this is the case, we will say that
$L$ is an RM pseudolattice. Denote by $O_K$ the ring of integers of
$K$. There exists a unique integer $f\ge 1$ (conductor) such that
$\roman{End}\,L =\bold{Z}+fO_K=:R_f,$ and $L$ is
a projective module of rank 1 over $R_f$. 
\smallskip
The module $L$ is endowed with a total ordering.

\smallskip

Every $K$, $f$
and a ordered projective module over $R_f$ come from a lattice.

\smallskip

(c${}^{\prime}$) If pseudolattices $L$ and $L^{\prime}$ have the same
$K$ and $f$, they are isomorphic if and only if
their classes in the Picard group  $\roman{Pic}\,R_f$ coincide. 

\smallskip

Unlike the case of lattices, the automorphism group
of a pseudolattice is always infinite, it is isomorphic
to $\bold{Z}\times \bold{Z}_2.$

\medskip

{\bf 3.5. Quantum tori as ``limits'' of elliptic curves.}
Comparison of the relevant geometric categories
suggests that two--dimensional quantum tori
can be thus considered as limits of elliptic curves.
More specifically, take a family
of Jacobi parametrized curves $E_{\tau}=\bold{C}/(e^{2\pi i\tau})$
with $\roman{Im}\,\tau > 0$ and $\tau \to \theta \in \bold{R}$.
It is then natural to imagine $T_{\theta}$ as a limit
of $E_{\tau}.$ 

\smallskip

Fixing a Jacobi uniformization of an elliptic curve
(or abelian variety of any dimension) as 
a part of its structure is necessary, for example, 
in problems connected with mirror symmetry.
In such contexts our intuition seemingly provides
a sound picture (cf. a similar discussion in [So],
pp. 100, 113--114).

\smallskip

However, limitations
of this viewpoint become quite apparent if one
has no reason to keep a Jacobi uniformization
as a part of the structure, and is interested only
in the isomorphim classes of elliptic curves, perhaps
somewhat rigidified by a choice of a level structure.

\smallskip

In this case one must contemplate the dynamics of the limiting process
not on the closed upper half--plane but on a relevant modular
curve $X$. Letting $\tau$ tend to $\theta$ along a geodesic,
we get a parametrized real curve on $X$ which, when $\theta$
is irrational, does not tend to any limiting point.
This is what can happen.

\medskip

(a) Let $\theta$ be a real quadratic
irrationality, $\theta^{\prime}$ its conjugate.
Consider the oriented geodesic in $H$ joining $\theta^{\prime}$
to $\theta$. The image of this geodesic on any modular curve
$X$ is supported by
a closed loop, which we denote $(\theta^{\prime},\theta )_X$.

\smallskip

(b) Let $\theta$ be as above, and let $\tau$ tend to $\theta$
along an arbitrary geodesic. Then the image of this
geodesic on $X$ has $(\theta^{\prime},\theta )_X$
as a limit cycle (in positive time).

\smallskip

(c) Each closed geodesic on $X$ is the support 
of a closed loop $(\theta^{\prime},\theta )_X$. The union
of them is dense in $X$. It is a strange attractor  
for the geodesic flow in the following sense.
Having chosen a sequence of loops $(\theta^{\prime}_i,\theta_i )_X$,
a sequence of integers $n_i\ge 1$, and a sequence of real
numbers $\epsilon_i>0$, $i=1,2, \dots $, one can find
an oriented geodesic winding $\ge n_i$ times in the 
$\epsilon_i$--neighborhood of $(\theta^{\prime}_i,\theta_i )_X$ for each
$i$, before jumping to the next loop.

\smallskip

Now let us imagine that we have constructed a certain object
$R(E_\tau )$ depending on the isomorphism class of $E_\tau$
(perhaps, with rigidity). This object can be a number,
a function of the lattice, a linear space, a category ... Suppose also
that we have constructed a similar object $\Cal{R} (T_{\theta})$
depending on the isomorphism class of $T_{\theta}$,
and that we want to make sense of the intuitive notion
that $\Cal{R} (T_{\theta})$ is ``a limit of $R(E_\tau )$.''
Since in the most interesting for us case (a) 
$E_{\tau}$ keeps rotating around the same loop,
there are  two natural possibilities:

\medskip

(i) {\it The object $R(E_\tau )$ actually ``does not depend on
$\tau$'', and $\Cal{R} (T_{\theta})$ is its constant
value. Here independence generally means a canonical identification
of different $R(E_\tau )$, e.g. via a version of  flat
connection defined along the loop.} 

\smallskip

(ii)  {\it The object $R(E_\tau )$ does depend on
$\tau$, and $\Cal{R} (T_{\theta})$ is obtained by a 
kind of integrating or averaging various
$R(E_\tau )$ along the loop.}

\medskip

The second case looks more interesting, however, it
is not immediately
obvious that such objects occur in nature.
Remarkably, they do, and precisely in the context
of real multiplication and Stark's conjecture.
In fact, this is how we will interpret the
beautiful old calculational tricks  due to Hecke:
see [He1], [He2], [Her], [Za1]. See also [Dar]
for a similar observation related to what Darmon
calls Stark--Heegner points of elliptic curves.

\smallskip

In this section we will
only explain the geometric meaning of Hecke's substitution,
whereas the (slightly generalized) calculation itself will be treated
in the next section. 

\smallskip

Let $K\subset \bold{R}$
be a real quadratic subfield of $\bold{R}$ and
$L\subset K$ an RM pseudolattice. From now on, we denote
by $l\mapsto l^{\prime}$ the nontrivial element
of the Galois group of $K/\bold{Q}$. 

\smallskip

For any real $t$, consider the following subset of $\bold{C}$:
$$
\Lambda_t = \Lambda_t (L) :=
\{\lambda_t=\lambda_t(l):=le^{t/2}+il^{\prime}e^{-t/2}\,|\, l\in L\}
$$

\smallskip

\proclaim{\quad 3.5.1. Lemma} (a) $\Lambda_t(L)$ is a lattice.

\smallskip

(b) Any isomorphism $a:\,L_1\to L$ in the narrow sense
induces  isomorphisms $\Lambda_t(L_1)\to \Lambda_{t+c}(L)$
where $c$ is a constant depending only on $a$
and $t$ is arbitrary.
\smallskip

(c) The image of the curve $\{\Lambda_t\,|\,t\in\bold{R}\}$
on any modular curve  is a closed
geodesic. The affine coordinate $t$ along this curve
is the geodesic length.
\endproclaim

\medskip

{\bf 3.5.2. Remark.} One can try to relate elliptic curves
to quantum tori by treating these curves themselves
as objects on noncommutative geometry represented
by some version of the relevant crossed product algebra.
In the most direct approach, the latter is a
completion of the non--unitary toric algebra generated
by $U,V$ with
$UV=e^{2\pi i\tau}VU.$ Representation of such an algebra
are in fact closely related to vector bundles
on $E_{\tau}$ as was shown in [BEG] (following [BG]).
Developing further this approach, one can hope to see better
what happens when one passes to the unitary limit
$\tau\to\theta$.

\bigskip

\centerline{\bf \S 4. Stark's numbers for real quadratic fields}

\medskip

In this section we will explain Stark's conjecture
for real quadratic fields and slightly generalize Hecke's
method of calculation of these numbers. Since it involves
integration along the geodesic loops introduced
above, we conclude that Stark's numbers
have something to do with the respective
quantum tori. However, key parts
of the picture are still missing.

\smallskip

{\bf 4.1. Stark's numbers at $s=0$.} In this section we fix
a real quadratic subfield $K\subset \bold{R}$. Denote
by $l\mapsto l^{\prime}$  the action of the nontrivial
element of the Galois group of $K$, and by $O_K$ the ring
of integers of $K$, and put $N(l)=ll^{\prime}.$

\smallskip

Let $L$ be an arbitrary integral ideal of $K$ which,
together with its embedding in $\bold{R}$ and the induced
ordering, will be considered as a pseudolattice.

\smallskip

Choose also an $l_0\in O_K$ so that the pair $(L,l_0)$
satisfies the following restrictions:

\medskip

(i) {\it The ideals $\frak{b}:=(L,l_0)$ and $\frak{a}_0:=(l_0)\frak{b}^{-1}$
are coprime with $\frak{f}:=L\frak{b}^{-1}$.}

\smallskip

(ii) {\it Let $\varepsilon$ be a unit of $K$ such that
$\varepsilon\equiv 1\,\roman{mod}\,\frak{f}.$ Then $\varepsilon^{\prime}>0.$}

\medskip

Put now
$$
\zeta (L,l_0,s):= \roman{sgn}\,l_0^{\prime}\,N(\frak{b})^s\sum^{(u)}_{l\in L}
\frac{\roman{sgn}\,(l_0+l)^{\prime}}{|N(l_0+l)|^s}
\eqno(4.1)
$$
where $(u)$ at the summation sign means that one should
take one representative from each coset $(l_0+l)\varepsilon$
where $\varepsilon$ runs over all units $\equiv\,1\,\roman{mod}\,\frak{f}$.
Notice that $(l_0+L)\varepsilon =l_0+L$ precisely for such units.

\smallskip

With this conventions, our $\zeta (L,l_0,s)$ is exactly
Stark's function denoted $\zeta(s,\frak{c})$
on the page 65 of [St1]: our $\frak{b},\frak{f}$
have the same meaning in [St1], and our $l_0$ is Stark's $\gamma$.
The meaning of Stark's $\frak{c}$ is explained below.

\smallskip

The Stark number of $(L,l_0)$ is defined as
$$
S_0(L,l_0):=e^{\zeta^{\prime}(L,l_0,0)}\, .
\eqno(4.2)
$$

\smallskip

The simplest examples correspond to the cases when 
$(L,l_0)=(1)$, $\frak{f}=L$, in particular, $l_0=1$.

\smallskip

Notice that pseudolattices which are integral ideals
have conductor $f=1$.

\medskip

{\bf 4.2. Stark's conjecture for real quadratic fields.} 
In [St1], Stark conjectures
that $S_0(L,l_0)$ are algebraic units generating
abelian extensions of $K$. To be more precise,
let us first describe an abelian extension $M/K$
associated with $(L,l_0)$ using the classical
language of class field theory. (Our $M$ is Stark's $K$,
whereas our $K$ corresponds to Stark's $k$.)

\smallskip

In 4.1 above we constructed, starting with $(L,l_0)$,
the ideals $\frak{f}$ and $\frak{b}$ in $O_K$.
Let $I(\frak{f})$ be the group of fractional ideals of $K$
generated by the prime ideals of $K$ not dividing $\frak{f}$,
and $S(\frak{f})$ be its subgroup called the
principal ray class modulo $\frak{f}$. Then Artin's 
reciprocity map identifies 
$G(\frak{f}):=I(\frak{f})/S(\frak{f})$ with the Galois group of $M/K$.

\smallskip

Consider all pairs $(L,l_0)$ as above with fixed $\frak{f}.$
It is not difficult to establish that on this set,
$S_0(L,l_0)$ in fact depends only on the class $\frak{c}$ of
$(l_0)\frak{b}^{-1}$ in $G(\frak{f})$. Denote the respective number $E(\frak{c}).$

\medskip

{\bf 4.2.1. Conjecture.} {\it The numbers $E(\frak{c})$ are units
belonging to $M$ and generating $M$ over $K$.
If the Artin isomorphism associates with $\frak{c}$ an
automorphism $\sigma$, we have $E(1)^{\sigma}=E(\frak{c}).$}

\smallskip

(We reproduced here the most optimistic form of the Conjecture 1
on page 65 of [St1] involving $m=1$ and Artin's  reciprocity map).

\medskip

{\bf 4.3. Hecke's formulas.}  In this section we will
work out Hecke's approach to the computation of sums
of the type (4.1), cf. [He2]. It starts with a
Mellin transform so that instead of Dirichlet series (4.1)
we will be dealing with a version of theta--functions
for real quadratic fields. We start with introducing
a class of such theta functions more general than
strictly needed for dealing with  (4.1) (and more general
than Hecke's one). 

\medskip

{\bf 4.4. Theta functions of pseudolattices.}
Let $K\subset \bold{R}$ be as in 4.1. We choose and fix the
following data: a pseudolattice $L\subset K$,
two numbers $l_0,m_0\in K$ and a number
$\eta =\eta_0+i\eta_1\in\bold{C}.$ A complex variable $v$ will take values in
the upper half plane; $\sqrt{-iv}$ is the branch which is positive
on the upper part of the imaginary axis.

\smallskip

Finally, choose an infinite cyclic group $U$ of totally
positive units
in $K$ such that the following conditions hold:

\smallskip

(a) $u(l_0+L)= l_0+L$ for all
$u\in U.$

\smallskip

(b) $\roman{tr}\,ulm_0\equiv\roman{tr}\,lm_0\,\roman{mod}\,\bold{Z}$,
$\roman{tr}\,ul_0m_0\equiv\roman{tr}\,l_0m_0\,\roman{mod}\,\bold{Z}$
for all $l\in L$, $u\in U,$ where $\roman{tr}:=\roman{tr}_{K/\bold{Q}}.$

\smallskip

 Let $\varepsilon >1$
be a generator of $U$.   

\smallskip

Put now
$$
\Theta_{L,\eta}^U
\thickfracwithdelims[]\thickness0{l_0}{m_0}(v):=
$$
$$
\sum_{l_0+l\,\roman{mod}\,U}
(\eta_0\,\roman{sgn}\,(l_0^{\prime}+l^{\prime}) +
\eta_1\,\roman{sgn}\,(l_0+l))\,
e^{2\pi i\,v|(l_0+l)(l_0^{\prime}+l^{\prime})|}
e^{-2\pi i\,\roman{tr}\,lm_0}
e^{-\pi i\,\roman{tr}\,l_0m_0} . 
\eqno(4.3)
$$
Notation $l_0+l\,\roman{mod}\,U$ means that we sum over
a system of representatives of orbits of $U$
acting upon $l_0+L$.

\smallskip

Notice that such $U$ always exists, and that if we choose
a smaller subgroup $V\subset U$, then
$$
\Theta_{L,\eta}^V
\thickfracwithdelims[]\thickness0{l_0}{m_0}(v)=
[U:V]\,\Theta_{L,\eta}^U
\thickfracwithdelims[]\thickness0{l_0}{m_0}(v).
$$

\smallskip

In order to relate these thetas to Stark's numbers,
consider the function
$$
\Theta_{L,1}^U
\thickfracwithdelims[]\thickness0{l_0}{0}(v)=
\sum_{l_0+l\,\roman{mod}\,U}
\roman{sgn}\,(l_0^{\prime}+l^{\prime})\,
e^{2\pi i\,v|(l_0+l)(l_0^{\prime}+l^{\prime})|} .
\eqno(4.4)
$$
Then we have
$$
\sum_{l_0+l\roman{mod}\,U}
\frac{\roman{sgn}\,(l_0^{\prime}+l^{\prime})}{|N(l_0+l)|^s}=
\frac{(2\pi)^s}{\Gamma (s)}\,
\int_0^{i\infty}(-iv)^{s}
\Theta_{L,1}^U
\thickfracwithdelims[]\thickness0{l_0}{0}(v)\,\frac{dv}{v} .
\eqno(4.5)
$$
We will now show that these RM thetas can be obtained
by averaging some theta constants (related to the complex lattices)
along the closed geodesics described above.

\medskip

{\bf 4.5. Theta constants along geodesics.} Starting with the same data
as in 4.4, we consider first of all the Hecke family
of lattices $\Lambda_t=\Lambda_t(L)$ (see 3.5 above).
From $l_0$ which was used  to shift $L$, we will produce
a shift of $\Lambda_t$:
$$
\lambda_{0,t}:=l_0\,e^{t/2}+il_0^{\prime}\,e^{-t/2}.
$$
The number $m_0$  determines a character of $L$
appearing in (3.3): $l\mapsto e^{-2\pi i\,\roman{tr}\,lm_0}$. 
Similarly, we will produce a character of $\Lambda_t$
from
$$
\mu_{0,t}:=m_0\,e^{t/2}+im_0^{\prime}\,e^{-t/2}
$$
by using the scalar product on $\bold{C}$
$$
(x\cdot y)=\roman{Im}\, xy= x_0y_1+x_1y_0 
\eqno(4.6)
$$
where $x=x_0+ix_1,\,y=y_0+iy_1.$ Since $l_0,m_0\in L\otimes\bold{Q}$,
we have similarly $\lambda_{0,t}, \mu_{0,t}\in\Lambda_t\otimes\bold{Q}$.
Omitting $t$ for brevity, we put:
$$
\theta_{\Lambda,\eta}
\thickfracwithdelims[]\thickness0{\lambda_0}{\mu_0}\,(v):=
\sum_{\lambda\in\Lambda}
((\lambda_{0}+\lambda)\cdot\eta )\,e^{\pi i v|\lambda_{0}+\lambda|^2}
e^{-2\pi i (\lambda\cdot \mu_{0})-\pi i (\lambda_{0}\cdot
\mu_{0})} .
\eqno(4.7)
$$
The two types of thetas are related by Hecke's averaging formula:

\medskip

\proclaim{\quad 4.6. Proposition} We have
$$
\Theta_{L,\eta}^U
\thickfracwithdelims[]\thickness0{l_0}{m_0}(v)=
\sqrt{-iv}\,
\int_{-\roman{log}\,\varepsilon}^{\roman{log}\,\varepsilon} 
\theta_{\Lambda_t,\eta}
\thickfracwithdelims[]\thickness0{\lambda_{0,t}}{\mu_{0,t}}\,(v)\, dt .
\eqno(4.8)
$$
\endproclaim

{\bf Proof.} The following formulas are valid
for $\roman{Im}\,v>0$:
$$
e^{2\pi i\,v|mm^{\prime}|} =
\sqrt{-iv}\,|m^{\prime}|\int_{-\infty}^{\infty}
e^{-t/2}e^{\pi iv(m^2e^{t}+ m^{\prime 2}e^{-t})}\,dt =
$$
$$
\sqrt{-iv}\,|m|\int_{-\infty}^{\infty}
e^{t/2}e^{\pi iv(m^2e^{t}+ m^{\prime 2}e^{-t})}\,dt
\eqno(4.9)
$$
(see e.~g. [La], pp. 270--271).  In the rhs of (4.3), 
replace the first exponent
by the  integral expressions (4.9), using the first version
at $\eta_0$ and the second at $\eta_1$. We get:
$$
\Theta_{L,\eta}^U
\thickfracwithdelims[]\thickness0{l_0}{m_0}(v)=
$$
$$
\sqrt{-iv}\int_{-\infty}^{\infty} 
\sum_{l_0+l\,\roman{mod}\,U}
(\eta_0\,(l_0^{\prime}+l^{\prime})\,e^{-t/2} +
\eta_1\,(l_0+l)\,e^{t/2})\,\times
$$
$$
e^{\pi iv((l_0+l)^2e^{t}+ (l^{\prime}_0+l^{\prime})^2e^{-t})}
e^{-2\pi i\roman{tr}\,lm_0}
e^{-\pi i\roman{tr}\,l_0m_0}\,dt\,.
\eqno(4.10)  
$$
In view of (4.6) we have
$$
\eta_0\,(l_0^{\prime}+l^{\prime})\,e^{-t/2} +
\eta_1\,(l_0+l)\,e^{t/2} =
((\lambda_{0,t}+\lambda_{t})\cdot\eta ),
$$
$$
(l_0+l)^2e^{t}+ (l^{\prime}_0+l^{\prime})^2e^{-t}=
|\lambda_{0,t}+\lambda_t|^2,
$$
and similarly
$$
\roman{tr}\,lm_0=(\lambda_t\cdot\mu_{0,t}),\quad
\roman{tr}\,l_{0}m_{0}=(\lambda_{0,t}\cdot\mu_{0,t}) .
$$
Inserting this into (4.10), we obtain
$$
\sqrt{-iv}\,\int_{-\infty}^{\infty} \,dt
\sum_{l_0+l\,\roman{mod}\,U}
((\lambda_{0,t}+\lambda_t)\cdot\eta)
e^{\pi i\, v|\lambda_{0,t}+\lambda_t|^2}\,
e^{-2\pi i\, (\lambda_t\cdot\mu_{0,t})}
e^{-\pi i\, (\lambda_{0,t}\cdot\mu_{0,t})} .
\eqno(4.11)  
$$
Replacing $l_0+l$ by $\varepsilon (l_0+l)$ is equivalent
to replacing $t$ by $t+2\,\roman{log}\,\varepsilon$.
Hence finally the right hand side of (4.11) can be rewritten
as
$$
\sqrt{-iv}\,
\int_{-\roman{log}\,\varepsilon}^{\roman{log}\,\varepsilon} dt 
\sum_{\lambda_t\in\Lambda_t}
((\lambda_{0,t}+\lambda_t)\cdot\eta )\,e^{\pi i\, v|\lambda_{0,t}+\lambda_t|^2}
e^{-2\pi i\, (\lambda_t\cdot \mu_{0,t})-\pi i\, (\lambda_{0,t}\cdot
\mu_{0,t})}
\eqno(4.12)
$$
which is the same as (4.8).

\smallskip

We will now apply Poisson formula in order to derive
functional equations for Hecke's thetas.

\medskip

{\bf 4.7. Poisson formula.} Let $V$ be a real vector space,
$\widehat{V}$ its dual. We will denote by $(x\cdot y)\in\bold{R}$ the scalar
product of $x\in V$ and $y\in \widehat{V}.$
Choose a lattice (discrete subgroup of finite covolume) $\Lambda\subset
V$ and put
$$
\Lambda^!:=\{\mu \in \widehat{V}\,|\,\forall \lambda \in \Lambda,\,
(\lambda\cdot\mu )\in\bold{Z}\} .
\eqno(4.13)
$$
Choose also a Haar measure $dx$ on $V$ and define the Fourier transform
of a Schwarz function $f$ on $V$ by
$$
\widehat{f}(y):=\int_V f(x)\,e^{-2\pi i (x\cdot y)} dx . 
\eqno(4.14)
$$
If $f(x)$ in this formula is replaced by
$f(x+x_0)\,e^{-2\pi i (x\cdot y_0)-\pi i (x_0\cdot y_0)}$
for some $x_0\in V,\,y_0\in \widehat{V}$,
its Fourier transform $\widehat{f}(y)$ gets replaced
by $\widehat{f}(y+y_0)\,e^{2\pi i (x_0\cdot y)+\pi i (x_0\cdot y_0)}.$

\smallskip

The Poisson formula reads
$$
\sum_{\lambda\in\Lambda} f(\lambda )=\frac{1}{\int_{V/\Lambda} dx}\,
\sum_{\mu\in\Lambda^!} \widehat{f}(\mu ) ,
\eqno(4.15)
$$
and for shifted functions as above
$$
\sum_{\lambda\in\Lambda} f(\lambda_0+\lambda)\,e^{-2\pi i (\lambda\cdot \mu_0)-\pi i (\lambda_0\cdot \mu_0)} =\frac{1}{\int_{V/\Lambda} dx}\,
\sum_{\mu\in\Lambda^!} \widehat{f}(\mu_0+\mu )\,e^{2\pi i (\lambda_0\cdot \mu )+\pi i (\lambda_0\cdot \mu_0) } .
\eqno(4.16)
$$
\medskip
{\bf 4.8. Functional equations for $\theta$ and $\Theta$.}
In order to transform (4.12) using the Poisson formula,
we put
$$
V=\bold{C}=\{x_0+ix_1\},\ \widehat{V}= \bold{C}=\{y_0+iy_1\},\
\eqno(4.17)
$$
and take (4.6) for the scalar product.

\smallskip

\proclaim{\quad 4.8.1. Lemma} Let the lattice $\Lambda_t\subset \bold{C}$
be the Hecke lattice. Then the dual lattice $\Lambda_t^!$ with respect
to the pairing (4.6) has the similar structure
$$
\Lambda_t^! = \Lambda_t (M) :=
\{me^{t/2}+im^{\prime}e^{-t/2}\,|\, m\in M\} 
\eqno(4.18)
$$
where we denoted by $M=L^?$ the pseudolattice
$$
M:=\{ m\in K\,|\,\forall l\in L,\,\roman{tr}_{K/\bold{Q}} (l^{\prime}m)
\in\bold{Z}.\}.
$$
\endproclaim
\smallskip

{\bf Proof.} Denote by $\Gamma$ the lattice (4.18).
For any $\lambda =  le^{t/2}+il^{\prime}e^{-t/2}\in 
\Lambda_t$ and $\mu = me^{t/2}+im^{\prime}e^{-t/2}\in \Gamma$
we have
$$
(\lambda\cdot\mu )= \roman{Im}\,\lambda\mu=lm^{\prime}+l^{\prime}m=
\roman{tr}_{K/\bold{Q}}(lm^{\prime}).
\eqno(4.19)
$$
Therefore this scalar product lies in $\bold{Z}$ if
$m\in M$ so that $\Gamma\subset \Lambda_t^!$.
Clearly, then, $\Gamma$ must be commensurable with $\Lambda_t^!$,
so that the right hand side of (4.19) can be used
for computing  $(\lambda\cdot\mu )$ on the whole $\Lambda_t^!$.
This finishes the proof.

\smallskip
For example,  $O_K^? =\frak{d}^{-1}$ where
$\frak{d}$ is the different. In fact,
this is the standard definition of the different. 

\smallskip

Now let $l_1,l_2$ be two generators of the pseudolattice $L$.
Put
$$
\Delta (L):= |l_1l_2^{\prime}-l_1^{\prime}l_2| .
\eqno(4.20)
$$
Clearly, this number does not depend on the choice of
generators.

\medskip

\proclaim{\quad 4.8.2. Lemma} Let the Haar measure on $V$ be
$dx=dx_0\,dx_1.$ Choose generators $l_1,l_2$ of $L$. Then
$$
\int_{V/\Lambda_t}dx = \Delta (L) .
\eqno(4.21)
$$
\endproclaim
\smallskip

{\bf Proof.} If $\Lambda_t$ is generated by $\omega_1,\omega_2$,
then the volume (4.21) equals 
$$
|\roman{Re}\,\omega_1 \,\roman{Im}\,\omega_2 -
\roman{Re}\,\omega_2 \,\roman{Im}\,\omega_1 |.
$$
Taking
$$
\omega_1=l_1e^{t/2}+il_1^{\prime}e^{-t/2},\quad
\omega_2=l_2e^{t/2}+il_2^{\prime}e^{-t/2},
$$
we get (4.20).

\medskip

\proclaim{\quad 4.8.3. Lemma} The Fourier transform of
$$
f_{v,\eta}(x):= (x\cdot\eta )\,e^{\pi i v|x|^2},\ \eta =\eta_0+i\eta_1
\eqno(4.22)
$$
equals
$$
{g}_{v,\eta}(y):=\frac{i}{v^2}\,(y\cdot i{\bar\eta})\,e^{-\frac{\pi i}{v}|y|^2}
\eqno(4.23)
$$
\endproclaim
\smallskip

{\bf Proof.} Putting $w=-iv$ we have
$$
f_{v,\eta}(x)=(x_0\eta_1+x_1\eta_0)\,e^{-\pi w(x_0^2+x_1^2)},
$$
so that its Fourier transform by (4.13) and (4.14) is
$$
\eta_1 \int_{-\infty}^{\infty} e^{-\pi wx_0^2}\,e^{-2\pi i x_0y_1}\, x_0\,dx_0\,\cdot\,
\int_{-\infty}^{\infty} e^{-\pi wx_1^2}\,e^{-2\pi i x_1y_0}\,dx_1 +
$$
$$
\eta_0 \int_{-\infty}^{\infty} e^{-\pi wx_0^2}\,e^{-2\pi i x_0y_1}\, dx_0\,\cdot\,
\int_{-\infty}^{\infty} e^{-\pi wx_1^2}\,e^{-2\pi i x_1y_0}\,x_1\,dx_1 =
$$
$$
(\eta_0y_0+\eta_1y_1)\,\frac{1}{iw^2}\,e^{-\pi\frac{y_0^2+y_1^2}{w}} .
$$
This is (4.23).

\medskip

{\bf 4.8.4. A functional equation for $\theta$.} Let us now 
write (4.16) for $f=f_{v,\eta}$ and $\Lambda_t$:
$$
\sum_{\lambda\in\Lambda_t}
((\lambda_{0,t}+\lambda)\cdot\eta )\,e^{\pi i v|\lambda_{0,t}+\lambda|^2}
e^{-2\pi i (\lambda\cdot \mu_{0,t})-\pi i (\lambda_{0,t}\cdot
\mu_{0,t})}=
$$
$$
\frac{i}{\Delta (L)\,v^2}\,\sum_{\mu\in\Lambda_t^!}
((\mu_{0,t}+\mu)\cdot i\bar{\eta} )\,e^{-\frac{\pi i} {v} |\mu_0+\mu |^2}
e^{2\pi i (\lambda_{0,t}\cdot \mu )+\pi i (\lambda_{0,t}\cdot
\mu_{0,t})} .
$$ 
In the notation (4.7) this means:
$$
\theta_{\Lambda_t,\eta}
\thickfracwithdelims[]\thickness0{\lambda_{0,t}}{\mu_{0,t}}\,(v)=
\frac{i}{\Delta (L)\,v^2}\,
\theta_{\Lambda_t^!,i\bar{\eta}}
\thickfracwithdelims[]\thickness0{\mu_{0,t}}{-\lambda_{0,t}}\,
\left(-\frac{1}{v}\right) .
\eqno(4.24)
$$
\smallskip

We now can establish a functional equation
for $\Theta^U$ as well:

\medskip

\proclaim{\quad 4.9. Proposition} We have
$$
\Theta_{L,\eta}^U
\thickfracwithdelims[]\thickness0{l_0}{m_0}(v)=
\frac{1}{\Delta (L)\,v}\,
\Theta_{L^?,i\bar{\eta}}^U
\thickfracwithdelims[]\thickness0{m_0}{-l_0}\left(-\frac{1}{v}\right) .
\eqno(4.25)
$$
\endproclaim
\smallskip

{\bf Proof.} This is a straightforward consequence of (4.8)
and (4.24).

\bigskip

\centerline{\bf References}

\medskip

[ArLR] J.~Arledge, M.~Laca, I.~Raeburn. {\it Crossed products by semigroups
of endomorphisms and the Toeplitz algebras of ordered groups.}
Documenta Math., 2 (1997), 115--138.

\smallskip

[BG] V.~Baranovsky, V.~Ginzburg. {\it Conjugacy classes
in loop groups and $G$--bundles on elliptic curves.}
Int. Math. Res. Notices, 15 (1996), 733--751.

\smallskip

[BEG] V.~Baranovsky, S.~Evens, V.~Ginzburg. {\it Representations
of quantum tori and double--affine Hecke algebras.} e--Print
math.RT/0005024

\smallskip

[Bo] F.~Boca. {\it Projections in rotation algebras and theta functions.}
Comm. Math. Phys., 202 (1999), 325--357.

\smallskip 

[BoCo] J.~Bost, A.~Connes. {\it Hecke algebras, type III factors and phase
transitions with spontaneous symmetry breaking in number theory.}
Selecta Math., 3 (1995), 411--457.

\smallskip

[BrGrRi] L.~Brown, P.~Green, M.~A.~Rieffel.
{\it Stable isomorphism and strong Morita equivalence
of $C^*$--algebras.} Pacific J.~Math., 71 (1977),
349--363.

\smallskip

[Coh1] P.~Cohen. {\it A $C^*$--dynamical system with Dedekind zeta partition function
and spontaneous symmetry breaking.}
J. de Th. de Nombres de Bordeaux, 11 (1999), 15--30.

\smallskip

[Coh2] P.~Cohen. {\it Quantum statistical mechanics and number theory.}
In: Algebraic Geometry: Hirzebruch 70. Contemp. Math.,
vol 241,  AMS, Pridence RA, 1999, 121--128.

\smallskip

[Co1] A.~Connes. {\it Noncommutative geometry.} Academic Press,
1994.

\smallskip

[Co2] A.~Connes. {\it $C^*$--algebras et g\'eom\'etrie
diff\'erentielle.} C.~R.~Ac.~Sci.~Paris, t.~290 (1980), 599--604.

\smallskip

[Co3] A.~Connes. {\it Trace formula in noncommutative geometry
and the zeros of the Riemann zeta function.} Selecta Math., New Ser., 5 (1999),
29--106.

\smallskip

[Co4] A.~Connes. {\it Noncommutative geometry and the Riemann
zeta function.} In: Mathematics: Frontiers and Perspectives,
ed. by V.~Arnold et al., AMS, 2000, 35--54.

\smallskip

[Co5] A.~Connes. {\it Noncommutative Geometry Year 2000.}
e--Print   math.QA/0011193

\smallskip

[CoDSch] A.~Connes, M.~Douglas, A.~Schwarz. {\it Noncommutative
geometry and Matrix theory: compactification on tori.}
Journ. of High Energy Physics, 2 (1998).

\smallskip

[Dar] H.~Darmon. {\it Stark--Heegner points over real quadratic fields.}
Contemp. Math., 210 (1998), 41--69.

\smallskip

[DE] {\it The Grothendieck Theory of Dessins d'Enfants
(Ed. by L.~Schneps).} London MS Lect. Note series, 200,
Cambridge University Press, 1994.

\smallskip
[Dr] V.~G.~Drinfeld. {\it On quasi--triangular quasi--Hopf
algebras and some groups closely associated
with $Gal (\overline{\bold{Q}}/\bold{Q})$.}
Algebra and Analysis 2:4 (1990); Leningrad Math. J.
2:4 (1991), 829--860.

\smallskip

[Ga] C.~F.~Gauss.  {\it Zur Kreistheilung.} In: Werke, Band 10, p.4.
Georg Olms Verlag, Hildsheim--New York, 1981.

\smallskip

[GoHaJo] F.~Goodman, P.~de la Harpe, V.~Jones. {\it Coxeter graphs and
towers of algebras.} Springer, 1989.

\smallskip

[HaL] D.~Harari, E.~Leichtnam. {\it Extension du
ph\'enom\`ene de brisure spontan\'ee de sym\'etrie
de Bost--Connes au cas des corps globaux
quelconques.} Selecta Math., 3 (1997), 205--243.

\smallskip

[He1] E.~Hecke. {\it \"Uber die Kroneckersche Grenzformel
f\"ur reelle quadratische K\"orper und die Klassenzahl
relativ--abelscher K\"orper.} Verhandl. d. 
Naturforschender Gesell. i. Basel, 28 (1917),
363--373. (Math. Werke, pp. 198--207, Vandenberg \& Ruprecht,
G\"ottingen, 1970).

\smallskip

[He2] E.~Hecke. {\it Zur Theorie der elliptischen  
Modulfunktionen.} Math.~Annalen, 97 (1926), 210--243.
(Math. Werke, pp. 428--460, Vandenberg \& Ruprecht,
G\"ottingen, 1970).

\smallskip

[Her] G.~Herglotz. {\it \"Uber die Kroneckersche Grenzformel
f\"ur reelle quadratische K\"orper I, II.} Berichte
\"uber die Verhandl. S\"achsischen Akad. der Wiss.
zu Leipzig, 75 (1923), 3--14, 31--37.

\smallskip

[Jo1] V.~Jones. {\it Index for subfactors.} Inv. Math., 72:1 (1983), 1--25.

\smallskip

[Jo2] V.~Jones. {\it Index for subrings of rings.}
Contemp. Math. 43, AMS(1985), 181--190.

\smallskip

[Jo3] V.~Jones. {\it Fusion en alg\`ebres de von Neumann et groupes de
lacets (d'apr\`es A.~Wassermann)}. S\'eminaire Bourbaki, no. 800
(Juin 1995), 20 pp.

\smallskip

[La] S.~Lang. {\it Elliptic Functions.} Addison--Wesley, 1973.

\smallskip

[LoS] P.~Lochak, L.~Schneps. {\it A cohomological interpretation
of the Grothendieck--Teichm\"uller group.} Inv. Math.,
127 (1997), 571--600.

\smallskip

[Ma1] Yu.~Manin. {\it  Quantized theta--functions.} In: Common
Trends in Mathematics and Quantum Field Theories (Kyoto, 1990), 
Progress of Theor. Phys. Supplement, 102 (1990), 219--228.

\smallskip

[Ma2] Yu.~Manin. {\it Mirror symmetry and quantization of abelian varieties.}
In: Moduli of Abelian Varieties, ed. by C.~Faber et al.,
Progress in Math., vol. 195, Birkh\"auser, 2001, 231--254.
e--Print  math.AG/0005143

\smallskip

[Ma3] Yu.~Manin. {\it Theta functions, quantum tori and Heisenberg groups}.
Lett. in Math. Physics, 56 (2001), 295--320.
e--Print math.AG/001119

\smallskip

[Ma4] Yu.~Manin. {\it Real multiplication and noncommutative
geometry.} e--Print math.AG/0202109.

\smallskip

[MaMar] Yu.~Manin, M.~Marcolli. {\it Continued fractions, modular symbols, and non-commutative geometry.} e--Print math.NT/0102006

\smallskip

[NiVa] D.~Nikshych, L.~Vainerman. {\it Finite quantum
groupoids and their applications.} e--Print math.QA/0006057.

\smallskip

[Ri1] M.~A.~Rieffel. {\it Strong Morita equivalence of
certain transformation group $C^*$--algebras.}
Math. Annalen, 222 (1976), 7--23.

\smallskip

[Ri2] M.~A.~Rieffel. {\it Von Neumann algebras associated with pairs of
lattices in Lie groups.} Math.~Ann., 257 (1981), 403--418.

\smallskip

[Ri3] M.~A.~Rieffel. {\it $C^*$--algebras associated with irrational rotations.}
Pacific J.~Math., 93 (1981), 415--429.

\smallskip

[Ri4] M.~A.~Rieffel. {\it The cancellation theorem for projective
modules over irrational rotation $C^*$--algebras.}
Proc.~Lond.~Math.~Soc. (3), 47 (1983), 285--303.

\smallskip

[Ri5] M.~A.~Rieffel. {\it Projective modules over higher--dimensional
non--commutative tori.} Can.~J.~Math., vol.~XL, No.~2 (1988), 257--338.

\smallskip

[Ri6] M.~A.~Rieffel. {\it Non--commutative tori --- a case
study of non--commutative differential manifolds.}
In: Cont.~Math., 105 (1990), 191--211.

\smallskip

[RiSch] M.~A.~Rieffel, A.~Schwarz. {\it Morita equivalence
of multidimensional non--commutative tori.} 
Int. J. Math., 10 (1999), 289--299. e--Print   math.QA/9803057

\smallskip

[Sch1] A.~Schwarz. {\it Morita equivalence and duality.}
Nucl.~Phys., B 534 (1998), 720--738.

\smallskip

[Sch2] A.~Schwarz. {\it Theta--functions on non--commutative tori.}
e--Print math/0107186

\smallskip

[Se] J.~P.~Serre. {\it Complex Multiplication.} In: 
Algebraic Number Fields, ed. by J.~Cassels, A.~Fr\"olich.
Academic Press, NY 1977, 293--296.

\smallskip

[So] Y.~Soibelman. {\it Quantum tori, mirror symmetry
and deformation theory.} Lett. in Math. Physics,
56 (2001), 99--125. e--Print math.QA/0011162.

\smallskip

[St1] H.~M.~Stark. {\it $L$--functions at $s=1$. III. Totally real
fields and Hilbert's Twelfth Problem.} Adv. Math., 22 (1976), 64--84.

\smallskip
[St2] H.~M.~Stark. {\it $L$--functions at $s=1$. IV. First
derivatives at $s=0$.} Adv. Math., 35 (1980), 197--235.

\smallskip

[Ste] P.~Stevenhagen. {\it Hilbert's 12th problem, Complex
Multiplication and Shimura reciprocity.} In: Class Field Theory --
Its Centenary and Prospect. Adv Studies in Pure Math.,
30 (2001), 161--176.

\smallskip

[Ta] J.~Tate. {\it Les conjectures de Stark sur les fonctions
$L$ d'Artin en $s=0$.} Progress in Math., vol. 47, Birkh\"auser, 1984.

\smallskip

[V] S.~Vladut. {\it Kronecker's Jugendtraum
and modular functions.} Gordon and Breach, 1991.

\smallskip

[Wa] Y.~Watatani. {\it Index for $C^*$--subalgebras.}
Mem. AMS, vol. 83, Nr. 424, Providence, RA, 1990.

\smallskip

[Wi] E.~Witten. {\it Overview of $K$--theory applied to strings.}
e--Print hep-th/0007175

\smallskip

[Za1] D.~Zagier. {\it A Kronecker limit formula for real quadratic field.}
Math. Ann., 213 (1975), 153--184.

\smallskip

[Za2] D.~Zagier. {\it Valeurs des fonctions zeta des corps
quadratiques r\'eels aux entiers n\'egatifs.} Ast\'erisque
41--42 (1977), 135--151.

\enddocument